\documentstyle[12pt,amssymb]{amsart} 
\def\SBIMSMark#1#2#3{
 \font\SBF=cmss10 at 10 true pt
 \font\SBI=cmssi10 at 10 true pt
 \setbox0=\hbox{\SBF Stony Brook IMS Preprint \##1}
 \setbox2=\hbox to \wd0{\hfil \SBI #2}
 \setbox4=\hbox to \wd0{\hfil \SBI #3}
 \setbox6=\hbox to \wd0{\hss
             \vbox{\hsize=\wd0 \parskip=0pt \baselineskip=10 true pt
                   \copy0 \break%
                   \copy2 \break%
                   \copy4 \break}}
 \dimen0=\ht6   \advance\dimen0 by \vsize \advance\dimen0 by 8 true pt
                \advance\dimen0 by -\pagetotal
 \dimen2=\hsize \advance\dimen2 by .25 true in
  \openin2=publishd.tex
  \ifeof2\setbox0=\hbox to 0pt{}
  \else 
     \setbox0=\hbox to 3.1 true in{
                \vbox to \ht6{\hsize=3 true in \parskip=0pt  \noindent  
                \input publishd.tex 
                \vfill}}
  \fi
  \closein2
  \ht0=0pt \dp0=0pt
 \ht6=0pt \dp6=0pt
 \setbox8=\vbox to \dimen0{\vfill \hbox to \dimen2{\copy0 \hss \copy6}}
 \ht8=0pt \dp8=0pt \wd8=0pt
 \copy8
 \message{*** Stony Brook IMS Preprint #1, #2 ***}
}

\renewcommand{\marginpar}[1]{}
\catcode`\@=12

\def\Empty{}
\newcommand\oplabel[1]{
  \def\OpArg{#1} \ifx \OpArg\Empty {} \else
  	\label{#1}
  \fi}
		
\long\def\realfig#1#2#3#4{
\begin{figure}[htp]
\centerline{\psfig{figure=#2,width=#4}}
\caption[#1]{#3}
\oplabel{#1}
\end{figure}}

\newcommand{\comm}[1]{}

\input{psfig} 
\newtheorem{thm}{Theorem}
\newtheorem{cor}{Corollary}
\newtheorem{lem}{Lemma}
\newtheorem{prop}{Proposition}

\newcommand{\thmref}[1]{Theorem~\ref{#1}}
\newcommand{\propref}[1]{Proposition~\ref{#1}}

\newcommand{\lemref}[1]{Lemma~\ref{#1}}
\newcommand{\corref}[1]{Corollary~\ref{#1}} 
\newcommand{\figref}[1]{Fig.~\ref{#1}}

\newcommand{\BBB}[1]{{\Bbb #1}}

\newcommand{\tinf}{{\BBB T}_{\infty}} 
 
\begin{document} 
\SBIMSMark{1998/1c}{January 1998}{}
\title{Biaccessiblility in Quadratic Julia Sets II: The Siegel and Cremer
Cases}\author{Saeed Zakeri}  
\address{Department of Mathematics, SUNY at Stony Brook, NY 11794} 
\email{zakeri@math.sunysb.edu} 
 
\pagestyle{myheadings} 
\markboth{\sc S. Zakeri}{\sc Biaccessiblility in Quadratic Julia Sets II } 
\begin{abstract} 
Let $f$ be a quadratic polynomial which has an irrationally indifferent
fixed point $\alpha$. Let $z$ be a biaccessible point in the Julia set of
$f$. Then: 
\begin{enumerate} 
\item[$\bullet$] 
In the Siegel case, the orbit of $z$ must eventually hit the critical point
of $f$. 
\item[$\bullet$] 
In the Cremer case, the orbit of $z$ must eventually hit the fixed point
$\alpha$. 
\end{enumerate} 
Siegel polynomials with biaccessible critical point certainly exist, but in
the Cremer case it is possible that biaccessible points can never exist.  
 
As a corollary, we conclude that the set of biaccessible points in the Julia
set of a Siegel or Cremer quadratic polynomial has Brolin measure zero.  
\end{abstract} 

\maketitle 

\noindent 
{\bf \S 1. Introduction.} Let $f$ be a polynomial map of the complex plane
$\BBB C$. A fixed point $z=f(z)$ is called {\it indifferent} if the {\it
multiplier} $\lambda= f'(z)$ has the form $e^{2\pi i \theta}$, where the
{\it rotation number} $\theta$ belongs to $\BBB R/ \BBB Z$. We call $z$ {\it
irrationally indifferent} if $\theta$ is irrational so that $\lambda$ is on
the unit circle but not a root of unity.  
 
Let $z$ be an irrationally indifferent fixed point of $f$. When $f$ is
holomorphically linearizable about $z$, we call $z$ a {\it Siegel} fixed
point. On the other hand, when $z$ is nonlinearizable, it is called a {\it
Cremer} fixed point.  
 
In this paper we only consider quadratic polynomials. Such a polynomial,
which we can put in the normal form  
\begin{equation} 
\label{eqn:newform} 
f: z\mapsto z^2+c, 
\end{equation} 
has two fixed points 
$(1\pm \sqrt{1-4c})/2$ 
which are distinct if and only if $c\neq 1/4$. If $c \notin [1/4, \infty)$,
so that the two fixed points have distinct real parts, then by convention
the fixed point which is further to the left is called $\alpha$ and the
other fixed point $1-\alpha$ is called $\beta$. The corresponding
multipliers are $\lambda=2 \alpha$ and $2-\lambda =2 \beta$, with
$|\lambda|< |2-\lambda|$. Evidently only the $\alpha$-fixed point can be
indifferent. The critical value parameter $c$ is then given by  
$$c=\lambda (2-\lambda)/4.$$   
Therefore, the set of all quadratic polynomials which have an indifferent
fixed point is a cardioid in the $c$-plane parametrized by $\lambda$ on the
unit circle. The set of quadratic polynomials with an irrationally
indifferent fixed point is then a dense subset of this cardioid.  
 
We call a quadratic polynomial $f$ in (\ref{eqn:newform}) {\it Siegel} or
{\it Cremer} if the $\alpha$-fixed point is irrationally indifferent and has
the corresponding property.     
 
By the theorem of Brjuno-Yoccoz \cite{Yoccoz}, $f$ is a Siegel polynomial if
and only if $\theta =\frac{1}{2 \pi i}\log f'(\alpha)$ satisfies the {\it
Brjuno condition}: 
\begin{equation} 
\label{eqn:Brjuno} 
\sum _{n=1}^{\infty} \frac{\log q_{n+1}}{q_n}< +\infty, 
\end{equation} 
where the $q_n$ appear as the denominators of the rational approximations
coming from the continued fraction expansion of $\theta$. 
 
Recall that the {\it filled Julia set} of $f$ is 
$$K(f)=\{ z\in \BBB C: \mbox{The orbit}\ \{ f^{\circ n}(z) \} _{n\geq 0}\
\mbox{is bounded} \}$$ 
and the {\it Julia set} of $f$ is the topological boundary of the filled
Julia set: 
$$J(f)=\partial K(f).$$ 
Both sets are nonempty, compact, connected and the filled Julia set is full,
i.e., the complement $\BBB C \backslash K(f)$ is connected. Every connected
component of the interior of $K(f)$ is a topological disk called a {\it
bounded Fatou component} of $f$. In the Siegel case, the component $S$ of
the interior of $K(f)$ containing the fixed point $\alpha$ is called the
{\it Siegel disk} of $f$ on which the action of $f$ is holomorphically
conjugate to the rigid rotation $z\mapsto e^{2 \pi i \theta}z$.  
 
Since $f(z)=f(-z)$ by (\ref{eqn:newform}), the Julia set $J(f)$ is invariant
under the $180^{\circ}$ rotation $\tau: z\mapsto -z$. If $U$ is an open
Jordan domain in the plane such that $\overline{U}\cap
\tau(\overline{U})=\emptyset$, it follows that $f$ is univalent in some
Jordan domain $V$ containing the closure $\overline U$. 
 
According to Sullivan, every bounded Fatou component must eventually map to
the immediate basin of attraction of an attracting periodic point, or to an
attracting petal for a parabolic periodic point, or to a periodic Siegel
disk for $f$ (see for example \cite{Milnor1}). On the other hand, by
\cite{Douady1} a polynomial of degree $d\geq 2$ can have at most $d-1$
nonrepelling periodic orbits. It follows that in the Siegel case, every
bounded Fatou component eventually maps to the Siegel disk $S$ centered at
$\alpha$. In the Cremer case, however, we simply conclude that $K(f)$ has no
interior, so that $K(f)=J(f)$.\\ \\  
{\bf \S 2. Accessibility.} Given a quadratic polynomial $f$ as in
(\ref{eqn:newform}) with connected Julia set, there exists a unique
conformal isomorphism 
$$\varphi: \overline{\BBB C} \backslash K(f)\rightarrow \overline{\BBB C}
\backslash \overline{\BBB D},$$ 
called the {\it B\"{o}ttcher map}, with $\varphi(\infty)=\infty$ and
$\varphi'(\infty)>0$, which conjugates $f$ to the squaring map: 
\begin{equation} 
\label{eqn:Bottcher} 
\varphi(f(z))=(\varphi(z))^2. 
\end{equation} 
The $\varphi$-preimages of the radial lines and circles centered at the
origin are called the {\it external rays} and {\it equipotentials} of
$K(f)$, respectively. The external ray $R_t$, by definition, is  
$$\varphi^{-1} \{ re^{2 \pi i t}: r>1 \} ,$$ 
where $t\in \BBB R / \BBB Z$ is called the {\it angle} of the ray. From
(\ref{eqn:Bottcher}) it follows that  
$$f(R_t)=R_{2t\ (mod\ 1)}.$$ 
 
More generally, let us consider an arbitrary compact, connected, full set
$K\subset \BBB C$, and let $\varphi_K: \overline{\BBB C} \backslash K
\rightarrow \overline{\BBB C} \backslash \overline{\BBB D}$ be the unique
conformal isomorphism with $\varphi_K(\infty)=\infty$ and
$\varphi'_K(\infty)>0$ given by the Riemann Mapping Theorem. We can define
the external rays $R_t$ for $K$ in a similar way. We say that $R_t$ {\it
lands} at $p\in \partial K$ if $\lim_{r\rightarrow 1}\varphi_K^{-1} (re^{2
\pi i t})=p$. A point $p\in \partial K$ is called {\it accessible} if there
exists a simple arc in $\BBB C \backslash K$ which starts at infinity and
terminates at $p$. According to a theorem of Lindel\"{o}f (see for example
\cite{Rudin}, p. 259), $p$ is accessible exactly when there exists an
external ray landing at $p$. We call $p$ {\it biaccessible} if it is
accessible through at least two distinct external rays. By a theorem of
F. and M. Riesz \cite{Milnor1}, $K \backslash \{ p \} $ is disconnected
whenever $p$ is biaccessible. It is interesting that the converse is also
true. More precisely, if there are at least $n>1$ connected components of $K
\backslash \{ p \} $, then at least $n$ distinct external rays land at $p$
(see for example \cite{McMullen}, p. 85). 
 
In the case $K=K(f)$ is the filled Julia set of a quadratic polynomial $f$
in (\ref{eqn:newform}), note that $\tau(R_t)=R_{t+1/2}$. Hence if $R_t$
lands at some $p\in J(f)$, then $R_{t+1/2}$ lands at $\tau(p)=-p$.\\ \\ 
{\bf \S 3. Arithmetical conditions.} As it is suggested by the theorem of
Brjuno-Yoccoz, the behavior of the orbits near the indifferent fixed point
is intimately connected to the arithmetical properties of the rotation
number $\theta$. There are certain classes of irrational numbers which are
of special interest in holomorphic dynamics and in this paper we will be
working with some of them. Let  
$$\theta=\frac{1}{a_1+\displaystyle{\frac{1}{a_2+\displaystyle{\frac{1}{\ddots}}}}}$$

be the continued fraction expansion of $\theta$, where all the $a_i$ are
positive integers, and  
$$\frac{p_n}{q_n}=\frac{1}{a_1+\displaystyle{\frac{1}{a_2+\displaystyle{\frac{1}{\ddots+\displaystyle{\frac{1}{a_n}}}}}}}$$

be the $n$-th rational approximation of $\theta$. We say that 
\begin{enumerate} 
\item [$\bullet$] 
$\theta$ is of {\it constant type} (we write $\theta \in {\cal {CT}}$) if
$\sup_n a_n < +\infty $.  
\item [$\bullet$] 
$\theta$ is {\it Diophantine} (we write $\theta \in \cal D$) if there 
exist positive constants $C$ and $\nu$ such that for every rational number
$0\leq p/q<1$, we have $|\theta -p/q|>C/q^{\nu}$. This condition is
equivalent to  
$\sup_n (\log q_{n+1}/\log q_n) < +\infty$. 
\item [$\bullet$] 
$\theta$ is of {\it Yoccoz type} (we write $\theta \in {\cal H}$) if every  
analytic circle diffeomorphism with rotation number $\theta$ is analytically

linearizable. (An explicit arithmetical description of $\cal H$ is given by
Yoccoz although it is not easy to explain; see \cite{Yoccoz}.) 
  
A closely related condition, which we denote by $\cal H'$, is defined as
follows: $\theta \in \cal H'$ if and only if every analytic circle
diffeomorphism with rotation number $\theta$, with no periodic orbit in some
neighborhood of the circle, is analytically linearizable
\cite{Perez-Marco1}. 
\item [$\bullet$] 
$\theta$ is of {\it Brjuno type} (we write $\theta \in {\cal B}$) if
$\theta$ satisfies the condition (\ref{eqn:Brjuno}). 
\end{enumerate} 
\vspace{0.1 in} 
We have the proper inclusions ${\cal H}\subset {\cal H'}$ and $\cal
{CT}\subset {\cal D} \subset {\cal H} \subset {\cal B}$. It is not hard to
show that $\cal D$, hence $\cal H, \cal H'$ and $\cal B$, are sets of full
measure in $\BBB R / \BBB Z$.\\ \\  
{\bf \S 4. Basic results.} Very little is known about the topology of the
Julia set of $f$ in the Siegel or Cremer case or the dynamics of $f$ on its
Julia set. The following theorem summarizes the basic results in the Cremer
case: 
\begin{thm} 
\label{Cremer} 
Let $f$ in $($\ref{eqn:newform}$)$ be a Cremer quadratic polynomial, so that
$\theta \notin \cal B$. Then 
\begin{enumerate} 
\item[(a)] 
The Julia set $J(f)$ cannot be locally-connected \cite{Sullivan}. 
\item[(b)] 
Every neighborhood of the Cremer fixed point $\alpha$ contains infinitely
many  
repelling periodic orbits of $f$ \cite{Yoccoz}. 
\item[(c)] 
The critical point $0$ is recurrent, i.e., it belongs to the closure of its
orbit  $\{ f^{\circ n}(0) \} _{n > 0}$ \cite{Mane}. 
\item[(d)] 
The critical point $0$ is not accessible from $\BBB C \backslash J(f)$
\cite{Kiwi}. 
\end{enumerate} 
\end{thm} 
See also \cite{Sorensen} for the so-called ``Douady's nonlanding Theorem''
which partially explains why the Julia set of a generic Cremer quadratic
polynomial fails to be locally-connected.  
 
In the Siegel case, we know a little bit more, but still the situation is
far from being fully understood. 
\begin{thm} 
\label{Siegel} 
Let $f$ in $($\ref{eqn:newform}$)$ be a Siegel quadratic polynomial, so that
$\theta \in \cal B$. Let $S$ denote the Siegel disk of $f$. Then 
\begin{enumerate} 
\item[(a)] 
If $\theta \in \cal{CT}$, then the boundary $\partial S$ is a quasicircle
which contains the critical point $0$ \cite{Douady2}. The Julia set $J(f)$
is locally-connected and has measure zero \cite{Petersen}. 
\item[(b)] 
If $\theta \in \cal H$, then $0\in \partial S$ \cite{Herman1}.  
\item[(c)] 
For some rotation numbers $\theta \in \cal B \backslash \cal H$, the entire
orbit of $0$ is disjoint from $\partial S$ \cite{Herman2}. In this case,
$J(f)$ cannot be locally-connected \cite{Douady2}. 
\item[(d)] 
For every $\theta \in \cal B$, the critical point $0$ is recurrent. 
\end{enumerate} 
\end{thm} 
Part (b) was proved by Herman for $\theta \in \cal D$, but his proof works
equally well for $\theta \in \cal H$. We will include a very short proof for
the latter case in section {\bf \S 5}. The proof of part (d) goes as
follows: If $0\in \partial S$, then by classical Fatou-Julia theory, every
point in $\partial S$ is in the closure of the orbit of $0$ \cite{Milnor1},
and recurrence follows. If $0\notin \partial S$ and $0$ is not recurrent,
then by \cite{Mane} the invariant set $\partial S$ is expanding, i.e., there
is a constant $\lambda >1$ and a positive integer $k$ such that $|(f^{\circ
k})'(z)|>\lambda$ for all $z\in \partial S$. It follows that the same
inequality holds over some neighborhood $U$ of $\partial S$ with a slightly
smaller $\lambda_1 >1$. We may as well assume that $U\cap S$ is
invariant. Take a small disk $V \Subset U\cap S$. Since $f^{\circ k}|_{U\cap
S}$ is holomorphically conjugate to the rigid rotation $z\mapsto e^{2 \pi i
k \theta}z$, there exists a sequence $n_j\rightarrow \infty$ such that
$f^{\circ k^{n_j}}$ converges uniformly to the identity map on $V$ as
$j\rightarrow \infty$. But this is impossible since for all $z\in V$,
$|(f^{\circ k^{n_j}})'(z)|>\lambda_1^{n_j}\rightarrow \infty$.\\ 
 
Comparing the two theorems, we notice that the Cremer case and the Siegel
case with $0 \notin \partial S$ share many properties. This is a general
philosophy which is partially explained by the theory of ``hedgehogs''
introduced recently by Perez-Marco \cite{Perez-Marco1} (see section {\bf \S
5} below). 
 
Inspired by this similarity, one expects the following to be true:\\ \\ 
{\bf Conjecture.} {\it Let $f$ be a Siegel quadratic polynomial and $0\notin
\partial S$. Then} 
\begin{enumerate} 
\item[(i)] 
{\it Every neighborhood of $\partial S$ contains infinitely many repelling
periodic orbits of $f$.} 
\item[(ii)]  
{\it The critical point $0$ is not accessible from $\BBB C\backslash
K(f)$.}\\ 
\end{enumerate} 
By an argument similar to \cite{Kiwi}, one can show that (i) implies (ii)
(see also \propref{H'}).\\ \\ 
{\bf \S 5. Hedgehogs.} Let $f$ be a Siegel or Cremer quadratic polynomial as
in (\ref{eqn:newform}). Let $U$ be a simply connected domain with compact
closure which contains the closure of the Siegel disk $S$ in the
linearizable case, or the indifferent fixed point $\alpha$ in the
nonlinearizable case. Suppose that $f$ is univalent in a neighborhood of the
closure $\overline U$. Then there exists a set $H=H_U$ with the following
properties: 
\begin{enumerate} 
\item[(i)] 
$\alpha \in H\subset \overline U,$ 
\item[(ii)] 
$H$ is compact, connected and full, 
\item[(iii)] 
$\partial H\cap \partial U$ is nonempty, 
\item[(iv)] 
$\partial H \subset J(f),$ 
\item[(v)]   
$f(H)=H.$ 
\end{enumerate} 
Note that $H$ has nonempty interior if and only if $\alpha $ is
linearizable. In this case our assumption that $f$ is univalent on $U$
implies that the critical point is off the boundary of the Siegel
disk. Clearly $H\supset \overline{S}$. 
 
Such an $H$ is called a {\it hedgehog} for the restriction
$f|_U:U\rightarrow \BBB C$. See \figref{hedgehog}(a) for the Cremer case and
(b) for the Siegel case. (We would like to emphasize that the topology of a
hedgehog is infinitely more complicated than anything we can possibly
sketch!) The existence of such totally invariant sets is proved by
Perez-Marco \cite{Perez-Marco1}.  
 
Note that in the Siegel case, one can get totally invariant sets $H$ with
the above  properties (i)-(v) even if $\partial U$ intersects the closure
$\overline S$. But in this case the existence of $H$ is not hard to show
because we can simply take $H$ as $\overline S$ or a compact invariant piece
with analytic boundary inside the Siegel disk (see \figref{hedgehog}(c) and
(d)).  
\realfig{hedgehog}{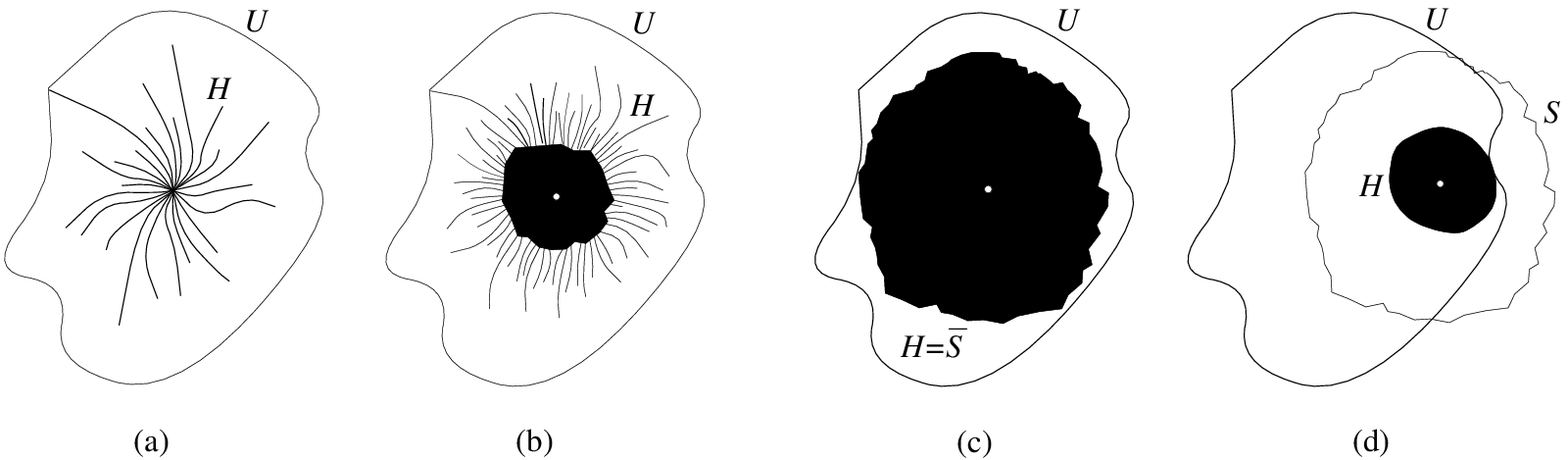}{}{13 cm} 
 
Hedgehogs turn out to be useful because of the following nice construction:
Uniformize the complement $\BBB C\backslash H$ by the Riemann map $\psi:\BBB
C \backslash H\rightarrow \BBB C \backslash \overline{\BBB D}$ and consider
the induced map $g=\psi\circ f \circ \psi^{-1}$ which is defined (by (v)
above) and holomorphic in an open annulus $\{ z\in \BBB C: 1< |z| < r
\}$. Use the Schwarz Reflection Principle to extend $g$ to the annulus $\{
z\in \BBB C: r^{-1}< |z| < r \}$. The restriction of $g$ to the unit circle
$\BBB T$ will then be a real-analytic diffeomorphism whose rotation number
is exactly $\frac{1}{2 \pi i}\log f'(\alpha )=\theta \in \BBB R/ \BBB Z$
(see \cite{Perez-Marco1}). This allows us to transfer results from the more
developed theory of circle diffeomorphisms to the less explored theory of
indifferent fixed points of holomorphic maps.  
 
Using the above construction, it is not hard to prove the following fact  
(see \cite{Perez-Marco2}): 
\begin{prop} 
\label{nobiacc} 
Let $p$ be a point in a hedgehog $H$ which is biaccessible from outside of
$H$. Then $p\in \partial S$ in the Siegel case and $p=\alpha$ in the Cremer
case. 
\end{prop} 
In fact, let us assume that we are in the Siegel case and $p\notin \partial
S$. Then one can find a simple arc $\gamma$ in $\BBB C\backslash H$ which
starts and terminates at $p$ and does not encircle the indifferent fixed
point $\alpha$. Let $D$ be the bounded connected component of $\BBB
C\backslash (H\cup \gamma)$. Evidently $\overline D$ is disjoint from
$\overline S$. The topological disk $D'=\psi(D)$ is bounded by the simple
arc $\psi(\gamma)$ and an interval $I$ on the unit circle. (The fact that
$\psi(\gamma)$ actually lands from both sides on the unit circle follows
from general theory of conformal mappings; see for example
\cite{Pommerenke}, page 29.) Since $g$ has irrational rotation number on the
unit circle $\BBB T$, for some integer $N$ we have $\bigcup_{i=0}^N g^{\circ
i}(I)=\BBB T$. By choosing $\gamma$ close enough to $H$, we can assume that
$g, g^{\circ 2}, \cdots, g^{\circ N}$ are all defined on $D'$ and
$\bigcup_{i=0}^N g^{\circ i}(D')$ contains an entire outer neighborhood of
$\BBB T$. It follows that $\bigcup_{i=0}^N f^{\circ i}(D)$ covers an entire
deleted neighborhood of $H$. Therefore, some iterate $f^{\circ i}(\overline
D)$ intersects $\partial S$. Since $f^{\circ i}$ is univalent on $\overline
D \cup \overline S$, it follows that $\overline D \cap \partial S \neq
\emptyset$, which contradicts our assumption. The proof in the Cremer case
is similar.\\  
 
The construction of the circle maps associated with hedgehogs as described
above gives short proofs for some interesting facts. As the first example,
we prove that there are no periodic points on $\partial S$ when the critical
point $0$ is off this boundary. This fact will be used in the proof of
\thmref{main} (see \cite{Perez-Marco1} for a general proof in the case of
indifferent germs; the fact that we are working with polynomials makes the
proof even shorter).  
 
First we need the following lemma: 
 
\begin{lem} 
\label{full} 
Let $f$ be a Siegel quadratic polynomial as in $($\ref{eqn:newform}$)$ whose
critical point $0$ is off the boundary $\partial S$ of the Siegel disk. Then
the closure $\overline S$ is full and $f$ acts injectively on it. 
\end{lem} 
 
It is reasonable to speculate that the closure of any bounded Fatou
component for a quadratic polynomial is full. This is known to be true
except when the polynomial has a periodic Siegel disk $S$ with the critical
point on its boundary $\partial S$. In this case, we do not know if
$\partial S$ can separate the plane into more than two connected components
(a so-called ``Lakes of Wada'' example in plane topology \cite{Hocking}).

\begin{pf} 
Since $f(z)=f(-z)$ for all $z$, if $f$ is not injective on $\overline S$,
there must be a pair of symmetric points $p$ and $-p=\tau(p)$ in $\partial
S$. Since $J(f)$ has a $180^{\circ}$ rotational symmetry, $f^{-1}(S)=S\cup
\tau (S)$. 
So $p$ and $-p$ also belong to $\partial (\tau(S))$. Consider the connected
component $V$ of $\BBB C \backslash (\overline{S} \cup \overline{\tau(S)})$
which contains the critical point $0$. Since $V$ is open and $\partial V
\subset J(f)$, it follows from the Maximum Principle that $V$ has to be a
bounded Fatou component of $f$. This contradicts the fact that $0\in J(f)$. 
  
Let us now assume that $\overline S$ is not full and let $U$ be a bounded
component of $\BBB C\backslash \overline S$. Since $\partial U \subset
\partial S \subset J(f)$, it follows again from the Maximum Principle that
$U$ has to be a bounded Fatou component of $f$, hence it eventually maps to
$S$, i.e., $f^{\circ n}(U)=S$ for some $n\geq 1$. Therefore $f^{\circ
n-1}(U)=\tau(S)$. But the boundary of $f^{\circ n-1}(U)$ is a subset of
$\partial S$, which implies that the common boundary $\partial S \cap
\partial (\tau (S))$ is nonempty. This contradicts the fact that
$f|_{\partial S}$ is injective. 
\end{pf} 
 
\begin{prop} 
\label{noperiod} 
Let $f$ be a Siegel quadratic polynomial whose critical point $0$ is off the
boundary $\partial S$. Then there are no periodic points on $\partial S$. 
\end{prop} 
 
\begin{pf} 
By the above lemma $\overline S$ is full and $f$ acts injectively on it, so
we can find a Jordan domain $U$ containing $\overline S$ such that
$f|_{\overline U}$ is univalent. Consider a hedgehog $H=H_U$ for the
restriction $f|_U$. Clearly $H\supset \overline S$. Suppose that there is a
periodic point on $\partial S$ which is necessarily repelling. Then there
exists a rational external ray $R$ landing at this point, hence $f^{\circ
n}(R)=R$ for some $n\geq 1$ (see for example \cite{Milnor1}). Consider the
induced map $g=\psi \circ f \circ \psi^{-1}$ as described above, and look at
the arc 
$\gamma= \psi (R)$. It is a standard fact that $\gamma$ has to land at some
point $p\in \BBB T$ \cite{Pommerenke} and $g^{\circ n}(p)=p$. But this
contradicts the fact that the rotation number of $g$ is irrational. 
\end{pf} 
 
In the second application, we prove \thmref{Siegel}(b): We want to show that
$\theta \in \cal H$ implies $0\in \partial S$. If not, by \lemref{full}
$\overline S$ is full and $f|_{\overline S}$ is univalent. Consider a Jordan
domain $U$, a hedgehog $H_U$ and the induced circle map $g$ as in the above
proof. Since the rotation number of $g$ belongs to $\cal H$, $g$ is
analytically linearizable. The linearization can be extended holomorphically
to an annulus neighborhood of the unit circle $\BBB T$. Pulling this
neighborhood back, we find a larger domain containing $S$ on which $f$ is
linearizable, which contradicts the definition of a Siegel disk. 
 
As a final application, we prove the following: 
 
\begin{prop} 
\label{H'} 
Let $f$ be a Siegel quadratic polynomial whose critical point $0$ is off the
boundary $\partial S$. If $\theta \in \cal{H}'$, then every neighborhood of
$\partial S$ contains infinitely many repelling periodic orbits. Hence, the
critical point $0$ is not accessible from $\BBB C\backslash K(f)$. 
\end{prop} 
 
\begin{pf} 
Consider the hedgehog construction as in the proof of \propref{noperiod} or
the above proof for \thmref{Siegel}(b). If there are no periodic orbits in
some neighborhood of $\partial S$, it follows that $g$ has no periodic orbit
in some neighborhood of $\BBB T$ either. Since the rotation number of $g$ is
$\theta \in \cal{H}'$, $g$ has to be linearizable. Now we can get a
contradiction as in the above proof for \thmref{Siegel}(b). So every
neighborhood of $\partial S$ must contain infinitely many periodic
orbits. The fact that this implies nonaccessibility of $0$ follows easily by
an argument similar to \cite{Kiwi}. 
\end{pf}   
\vspace{0.17in}   
{\bf \S 6. Wakes.} To see the behavior of rays near infinity, it will be
convenient to add a circle at infinity $\tinf \simeq \BBB R/\BBB Z$ to the
complex plane to obtain a closed disk $\copyright$ topologized in the
natural way. We denote the point $\lim_{r\rightarrow \infty} re^{2 \pi i t}$
on $\tinf$ simply by $\infty \cdot e^{2 \pi i t}$. The action of $f$ in
(\ref{eqn:newform}) on the complex plane extends continuously to
$\copyright$ by 
\begin{equation} 
\label{eqn:doubling}  
f(\infty \cdot e^{2 \pi i t})=\infty \cdot e^{4 \pi i t}, 
\end{equation} 
which is just the doubling map on $\tinf$. Note that the symmetry
$f(z)=f(-z)$ also extends to $\copyright$ if we define $- \infty \cdot e^{2
\pi i t}=\infty \cdot e^{2 \pi i (t+1/2)}$.\\ \\ 
{\bf Definition.} Let $f$ be a quadratic polynomial as in
(\ref{eqn:newform}) with connected Julia set. Let $z\neq \alpha$ be a
biaccessible point in $J(f)$ with two distinct rays $R$ and $R'$ landing on
it. We call $(R,R')$ a {\it ray pair}. By the Jordan Curve Theorem, $R \cup
R' \cup \{ z \}$ cuts the plane into two open topological disks. By the {\it
wake} $W$ of the ray pair $(R,R')$ we mean the connected component of $\BBB
C\backslash ( R \cup R' \cup \{ z \} )$ which does not contain the fixed
point $\alpha$. The other component is called the {\it co-wake} and it is
denoted by $\check{W}$. Point $z$ is called the {\it root} of $W$. The {\it
angle} $a(W)$ of the wake is just the (normalized) measure of
$\overline{W}\cap \tinf$. Clearly $a(W)+a(\check{W})=1$ (see
\figref{wake}(a)). 
 
Since distinct external rays are disjoint, it follows that any two wakes
with distinct roots are either disjoint or nested. 
 
In the following lemma we collect basic properties of wakes (compare with
\cite{Goldberg-Milnor} or \cite{Milnor2}): 
\begin{lem} 
\label{wakes} 
Let $z\in J(f)$ be a biaccessible point, $z\neq \alpha$, and let $W$ be a
wake with root $z$. 
\begin{enumerate} 
\item[(a)] 
If $z\neq 0$, then $a(W)>1/2$ if and only if $W$ contains the critical point
$0$.  
\item[(b)] 
If $a(W)=1/2$, then $z$ must be the critical point $0$. Conversely, if there
is any ray $R$ landing at $0$, then $R'=\tau(R)$ also lands at $0$ and the
two rays span a wake $W$ with $a(W)=1/2$.  
\item[(c)] 
Let $a(W)<1/2$ and $f(z)\neq \alpha$. Then $f(W)$ is a wake or a co-wake
with root $f(z)$, depending on whether $-\alpha \notin W$ or $-\alpha \in
W$. In either case, $f:W\rightarrow f(W)$ is a conformal  
isomorphism and $a(f(W))=2a(W)$.  
\end{enumerate} 
\end{lem} 
\begin{pf} 
Let $W$ be the wake of a ray pair $(R,R')$. 
  
(a) Let $0\in W$ and $a(W)<1/2$. Consider the symmetric region $\tau (W)$
whose angle is equal to $a(W)$. $W$ and $\tau (W)$ intersect since both
contain $0$ (see \figref{wake}(b)). On the other hand, $\overline{W} \cap
\overline{\tau(W)} \cap \tinf = \emptyset$ because $a(W)<1/2$. Since
$\overline{W}$ and $\overline{\tau(W)}$ are both homeomorphic to closed
disks, it follows that the ray pairs $(R,R')$ and $(\tau(R), \tau(R'))$ must
intersect, which is a contradiction. Therefore $a(W)>1/2$ if $0\in W$. 
 
On the other hand, let $a(W)>1/2$. Then the angle of the co-wake $\check{W}$
has to be less than $1/2$, so by the above argument $0 \notin \check{W}$, or
$0\in W$. This proves (a). 
 
\realfig{wake}{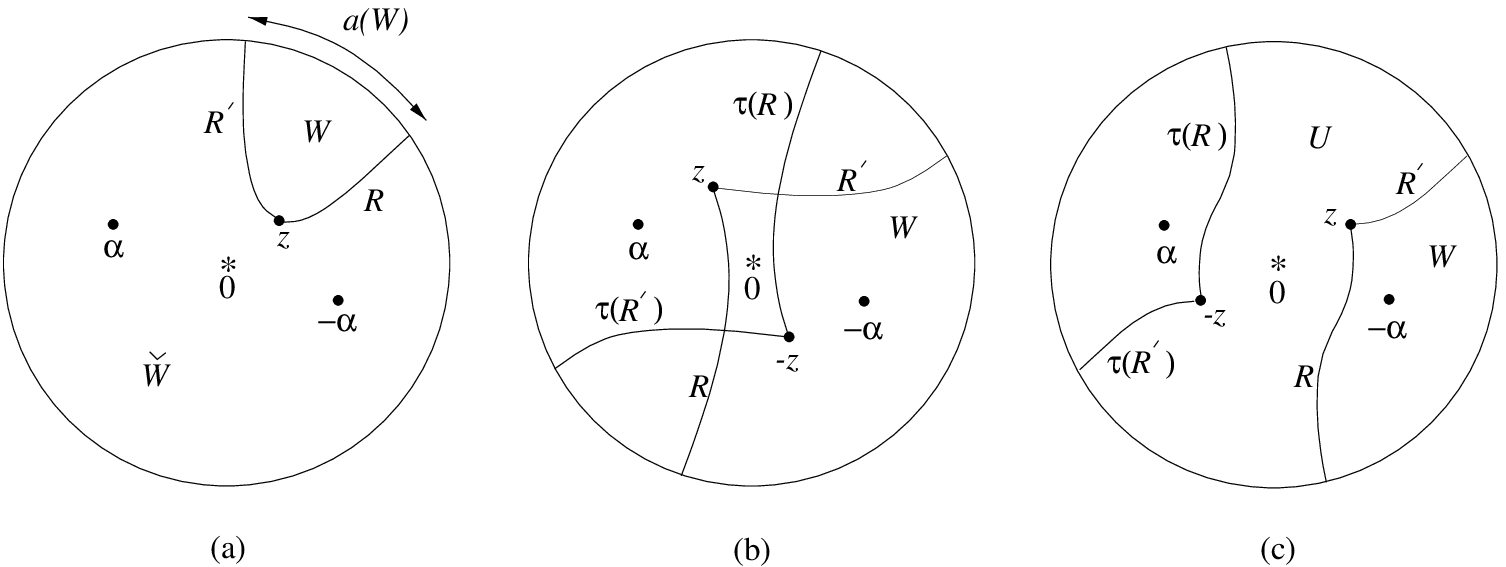}{}{15 cm}  
 
(b) If $a(W)=1/2$, then $R'=\tau(R)$. Hence $z=\tau(z)$ by continuity, which
means $z=0$. The converse is trivial. 
  
(c) If $a(W)<1/2$, then the ray pairs $(R,R')$ and $(\tau(R), \tau(R'))$ cut
the plane into simply connected domains $W$, $\tau(W)$ and an open set $U$
which is either a simply connected domain or the disjoint union of two
simply connected domains depending on whether $z\neq 0$ or $z=0$. By (a),  
$0\notin W \cup \tau(W)$. Consider the ray pair $(f(R), f(R'))$ landing at
$f(z)$, and let $W'$ be the corresponding wake. The pull-back of $W'$ by $f$
either consists of the disjoint union $W\sqcup \tau(W)$ or the open set $U$
(see \figref{wake}(c)). In the first case, $f$ maps $W$ to $W'$
isomorphically and $-\alpha \notin W$. In the second case, however, we must
have $-\alpha \in W$, $\alpha \in \tau(W)$, and both $W$ and $\tau(W)$ map
isomorphically to the co-wake $\check{W}'$. The fact that $a(f(W))=2a(W)$
simply follows from (\ref{eqn:doubling}).   
\end{pf} 
\vspace{0.17 in} 
{\bf \S 7. The main theorem.} Now we are in a position to state and  
prove the main theorem of this paper: 
 
\begin{thm} 
\label{main} 
Let $f$ be a quadratic polynomial as in $($\ref{eqn:newform}$)$ which has an
irrationally indifferent fixed point $\alpha$. Let $z$ be a biaccessible
point in the Julia set of $f$. Then: 
\begin{enumerate} 
\item[$\bullet$] 
In the Siegel case, the orbit of $z$ must eventually hit the critical point
$0$. 
\item[$\bullet$] 
In the Cremer case, the orbit of $z$ must eventually hit the fixed point
$\alpha$. 
\end{enumerate} 
\end{thm} 
 
(Compare \cite{Schleicher-Zakeri} where this same result for the Cremer case
is proved by a somewhat different argument.) 
 
In the Siegel case, if the critical point $0$ is accessible, then exactly
two rays land there (see the proof of \lemref{wakes}(b)). This happens, for
example, when $\theta \in \cal{CT}$, since in this case by
\thmref{Siegel}(a) the Julia set is locally-connected. On the other hand,
for some rotation numbers $\theta \in \cal{B} \cap \cal{H}'$, the critical
point is not accessible so that there are no biaccessible points in the
Julia set (see \corref{no}). 
 
In the Cremer case, if the fixed point $\alpha$ is accessible, then
infinitely many rays land there. In fact, if $R_t$ lands at $\alpha$, then
$t$ is irrational and every $R_{2^n t}$ lands at $\alpha$ also. However,
there is no known example where one can decide whether $\alpha$ is
accessible or not. 
 
\thmref{main} can be viewed from a more general perspective: Let $f$ be any
quadratic polynomial with connected Julia set $J(f)$. There is a unique
measure $\mu$ of maximal entropy $\log 2$, called the {\it Brolin measure},
which is supported on $J(f)$. In fact, $\mu$ coincides with the harmonic
measure induced by the radial limits of the inverse of the B\"{o}ttcher map
(see for example \cite{Lyubich}). It is a standard fact that $\mu$ is
$f$-invariant and ergodic. For $z\in J(f)$, let $v(z)$ denote the number of
external rays which land at $z$. (In Milnor's terminology \cite{Milnor2},
this is called the {\it valence} of $z$.) For $0\leq n \leq \infty$ define
$J_n=\{ z\in J(f): v(z)=n \}$. It follows from elementary plane topology
that the union $\bigcup_{n\geq 3} J_n$ is at most countable (see
\cite{Pommerenke}, page 36). On the other hand, the fact that almost every
external ray (with respect to the Lebesgue measure on $\BBB R/ \BBB Z$)
lands shows that $\mu(J_0)=0$. Putting these two facts together, we conclude
that $J(f)=J_1\cup J_2$ up to a set of $\mu$-measure zero. Note that
$v(f(z))=v(z)$ unless $z$ is the critical point. Therefore, if we neglect
the grand orbit of the critical point which has $\mu$-measure zero, it
follows that both $J_1$ and $J_2$ must be $f$-invariant subsets of the Julia
set. Ergodicity of $\mu$ then shows that up to a set of $\mu$-measure zero,
either $J(f)=J_1$ or $J(f)=J_2$. For example, in the first part of this
paper \cite{Zakeri}, it is shown that for any quadratic polynomial $f$ with
locally-connected Julia set, $J(f)=J_1$ must be the case unless $f$ is the
{\it Chebyshev} map $z\mapsto z^2-2$. For this map the Julia set is the
closed interval $[-2,2]$ and every point is the landing point of exactly two
rays except for the endpoints $\pm 2$ where unique rays land, so that
$J(f)=J_2$. \thmref{main} in particular proves that if $f$ has an
irrationally indifferent fixed point, then $J=J_1$ up to a set of Brolin
measure zero, thus covering some non locally-connected cases. It is
conjectured that $J=J_1$ is true for {\it every} non-Chebyshev quadratic
polynomial.\\     
  
The proof of \thmref{main} is based on the following lemma: 
 
\begin{lem} 
\label{mainlem} 
Let $f$ be a Siegel or Cremer quadratic polynomial as in
$($\ref{eqn:newform}$)$. Assume that there exists a biaccessible point in
$J(f)$ whose orbit never hits the critical point $0$ or the fixed point
$\alpha$. Then there exists a ray pair which separates $\alpha$ from $0$. 
\end{lem} 
 
\begin{pf} 
Let $z\in J(f)$ be such a biaccessible point and $(R,R')$ be a ray pair
which lands at $z$. Consider the associated wake $W_0$ with root $z$. Since
$z\neq 0$, we have $a(W_0)\neq 1/2$ by \lemref{wakes}(b). If $a(W_0)>1/2$,
then $0\in W_0$ by \lemref{wakes}(a) and $(R,R')$ separates $\alpha$ from
$0$. Let us consider the case where $a(W_0)<1/2$. If $-\alpha \in W_0$, then
$(R,R')$ must separate $-\alpha$ from $0$ because by \lemref{wakes}(a), $0
\notin W_0$. It follows that the symmetric ray pair $(\tau(R), \tau(R'))$
separates $\alpha$ from $0$. If, however, $-\alpha \notin W_0$, then by
\lemref{wakes}(c), $W_1=f(W_0)$ is a wake with root $z_1=f(z)$ with angle
$a(W_1)=2a(W_0)$.  
 
Now we can replace $W_0$ by $W_1$ in the above argument. If either
$a(W_1)>1/2$ or $a(W_1)<1/2$ and $-\alpha \in W_1$, we can find a ray pair
separating $\alpha$ from $0$. Otherwise, we consider the new wake
$W_2=f(W_1)$ with angle $a(W_2)=2^2 a(W_0)$. Since each passage $W_i\mapsto
W_{i+1}$ implies doubling the angles, this process must stop at some stage,
and this proves the lemma. 
\end{pf} 
\vspace{0.1 in} 
\noindent 
{\it Proof of \thmref{main}.} 
It will be more convenient to consider the Cremer case first. Suppose that
the orbit of $z$ never hits $\alpha$. Since the critical point is not
accessible by \thmref{Cremer}(d), \lemref{mainlem} gives us a ray pair
$(R,R')$ landing at some point $p\in J(f)$ which separates $\alpha$ from
$0$. Let $W$ be the corresponding wake with root $p$ and consider the
co-wake $\check W$. The restriction of $f$ to the closure of $\check W$ is
univalent since otherwise this closure would intersect the closure of the
symmetric domain $\tau (\check{W})$, which is impossible since
$a(\check{W})<1/2$. To work with a Jordan domain in the plane we cut off
$\check W$ along an equipotential curve and call the resulting domain $U$
(see \figref{proof}(a)). 
\realfig{proof}{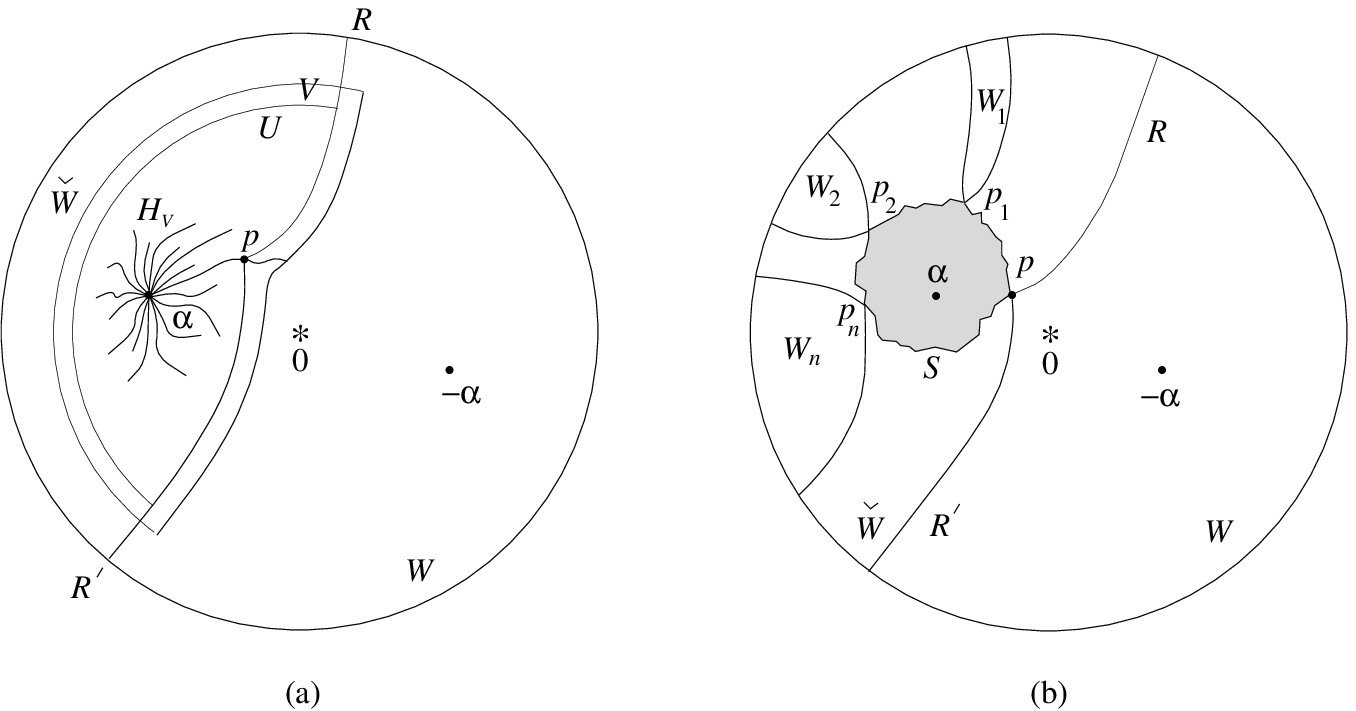}{}{13 cm}  
      
Now we consider a hedgehog $H_U$ for the restriction $f|_U: U\rightarrow
\BBB C$ as given in section {\bf \S 5}. Note that $p$ is the only point of
the Julia set on the boundary of $U$ and that $H_U\subset J(f)$ has to
intersect this boundary. Therefore, we simply have $p\in H_U$. Let us
consider a slightly larger Jordan domain $V \supset \overline U$ with
compact closure such that $f|_{\overline V}$ is still univalent. The
hedgehog $H_V$ for the restriction $f|_V: V\rightarrow \BBB C$ has to
contain $p$ also and reach the boundary of $V$. Since $p$ is biaccessible
from outside of the Julia set, it follows that $H_V \backslash \{ p \}$ is
disconnected. Therefore, $p$ is biaccessible from outside of $H_V$. This
contradicts \propref{nobiacc}, and finishes the proof of the theorem in the
Cremer case. 
 
Let us now assume that we are in the Siegel case. If the orbit of $z$
eventually hits the critical point $0$, there is nothing to
prove. Otherwise, since this orbit trivially cannot hit the fixed point
$\alpha \in S$, we are again in the situation of
\lemref{mainlem}. Therefore, there exists a ray pair $(R,R')$ landing at a
point $p\in J(f)$ which separates $\alpha$ from $0$. In particular the
critical point $0$ is off the boundary $\partial S$ of the Siegel disk. Then
the same argument as in the Cremer case with an application of
\propref{nobiacc} shows that $p$ must belong to $\partial S$. 
 
As before, let $W$ be the wake of the ray pair $(R,R')$, with root $p$. Then
by construction $W$ contains the critical point $0$ while the co-wake
$\check W$ contains the Siegel disk $S$ and has its boundary touching
$\overline S$ only at $p$. The point $p$ is not periodic by
\propref{noperiod}. Hence the successive images $p_n = f^{\circ n}(p) \in
\partial S$ are all contained in $\check W$ for $n \ge 1$. Therefore each
wake $W_n$ corresponding to the ray pair $(f^{\circ n}(R) , f^{\circ
n}(R'))$, with root point $p_n$, is also contained in $\check W$ (see
\figref{proof}(b)). In particular, none of these wakes contains the critical
point. Hence  $a(W_{n+1}) = 2 a(W_n) < 1/2$ for all $n$ by
\lemref{wakes}(c), which is clearly impossible. The contradiction shows that
the orbit of $z$ must eventually hit the critical point. \hfill$\Box$\\  
 
By \propref{H'}, we have the following corollary: 
 
\begin{cor} 
\label{no} 
Let $f$ be a Siegel quadratic polynomial with $0\notin \partial S$ and
$\theta \in \cal H'$. Then there are no biaccessible points in $J(f)$ at
all. 
\end{cor} 
 
By \lemref{wakes}(b), every wake with angle $1/2$ must have its root at the
critical point $0$. The converse is not true for arbitrary quadratic
polynomials. For example, the real Feigenbaum map $z\mapsto
z^2-1.401155\cdots$ has four distinct external rays landing on its critical
point (compare with \cite{Jiang}). 
However, in the case of a Siegel quadratic polynomial, the critical point
$0$ is the landing point of {\it at most} one ray pair $(R_t, \tau(R_t))$
(In the Cremer case, there are no such ray pairs by
\thmref{Cremer}(d)). This is nontrivial and follows from the statement that
every Siegel or Cremer quadratic on the boundary of the main cardioid of the
Mandelbrot set is the landing point of a unique parameter ray
\cite{Goldberg-Milnor}. In fact, one can explicitly compute the angle of the
candidate ray pair $(R_t, \tau(R_t))$ which may land at $0$ by 
$$t=\sum_{0 < p/q < \theta} 2^{-(q+1)}.$$ 
It is interesting that the uniqueness of such $t$ also follows from
\thmref{main}: 
 
\begin{cor} 
\label{uniquet} 
Let $f$ be a Siegel quadratic polynomial as in
$($\ref{eqn:newform}$)$. Then, no point in the Julia set $J(f)$ is the
landing point of more than two rays. In particular, at most one ray pair
lands at the critical point $0$. 
\end{cor} 
 
\begin{pf} 
By \thmref{main} it suffices to prove the corollary for the critical
point. Suppose that there is a ray pair $(R,R')$ which lands at $0$ such
that $R'\neq \tau (R)$. It follows that $(f(R), f(R'))$ is a ray pair which
lands at the critical value $c$. By \thmref{main}, the orbit of $c$ must
eventually hit the critical point $0$. But this means that $0$ is periodic,
which is impossible. 
\end{pf} 
\vspace{0.15 in}      
{\bf Acknowledgements.} This paper is an expanded version of a joint note
with Dierk Schleicher which was written during my visit to the University of
M\"{u}nich in July 1997, and covered the Cremer case with a rather different
proof \cite{Schleicher-Zakeri}. I would like to thank him for sharing his
ideas which led to the joint work, which have been useful in the present
paper as well. I am also grateful to Jack Milnor for his support and
interest in this work. In particular, he suggested the idea of working with
wakes, which simplified and unified part of the proof for both the Siegel
and Cremer cases.\\ \\

\end{document}